\documentclass[11pt]{amsart}
\catcode`\@=11 \catcode`\@=12

\usepackage{amsmath}
\usepackage{amssymb}
\usepackage{amsthm}
\usepackage{oldgerm}

\usepackage[pdftex,colorlinks,urlcolor=blue,pdfstartview=FitH]{hyperref}

\newcommand{\Rm}{\mathbb{R}}

\newcommand{\mF}{\ensuremath{\mathcal{F}}}

\newcommand{\mA}{\ensuremath{\mathcal{A}}}

\newcommand{\Tm}{\ensuremath{\mathbb{T}}}
\newcommand{\vs}{\vspace{.2cm}}

\newtheorem{lem}{Lemma}

\newtheorem{cor}[lem]{Corollary}
\newtheorem{prop}[lem]{Proposition}

\def\proof {\noindent{\sc{Proof. }}}
\def\qed {\mbox{}\hfill {\small \fbox{}} \\}
\def\lto{\longrightarrow}
\def\lmto{\longmapsto}

\def\leq{\leqslant}
\def\geq{\geqslant}

\title[Sub-solutions of the Hamilton-Jacobi equation]{Existence of $C^{1,1}$ critical 
sub-solutions of the  Hamilton-Jacobi equation
on compact manifolds}
\author{Patrick  Bernard}
\begin{document}

\maketitle
\begin{small}
\vs
 Patrick Bernard\\
CEREMADE, UMR CNRS 7534\\
Universit\'e de Paris Dauphine\\
Pl. du Mar\'echal de Lattre de Tassigny\\
75775 Paris Cedex 16,
France\\
\texttt{patrick.bernard@ceremade.dauphine.fr}\\
\end{small}
\vs
\hrule
\vs
R\'esum\'e: 
Nous donnons une preuve simple de l'existence 
d'une sous-solution $C^{1,1}$  de l'\'equation 
de Hamilton-Jacobi dans le contexte de la theorie
de Mather. Nous donnons  certaines cons\'equences
dynamiques de ce r\'esultat. Nous montrons que la solution
peut \^etre obtenue stricte en dehors de l'ensemble d'Aubry.

\vs
\hrule
\vs
Abstract: 
We offer a simple proof of the existence of a $C^{1,1}$ 
solution of the Hamilton-Jacobi equation in the context of Mather
theory. We derive some dynamical consequences of this result.
We also prove that the solution can be obtained strict outside 
of the Aubry set.\\
\vs
\hrule
\vs
Online access to this text:
\texttt{http://www.ceremade.dauphine.fr/ \~{}pbernard/}
\vspace{8cm}\\
This is a slightly corrected version of the text published under the same name in
Annales Scientifiques de l'\'Ecole Normale Sup\'erieure,
40, No.3, (2007), 445-452. Laurent Nocquet pointed out a  mistake in the proof 
of Proposition 9 in the published version, 
which is identical to Proposition 9 in the present version.
The present text, which solves that problem, has been posted on internet
in November 2010.

\newpage

Let $M$ be a compact manifold without boundary.
A function $H(x,p):T^*M\lto \Rm$ is called a Tonelli
Hamiltonian if it is $C^2$
and if, 
 for each $x\in  M$,
 the function $p\lmto H(x,p)$
 is convex with positive definite Hessian and superlinear
 on the fibre $T_x^*M$.
 Each Tonelli Hamiltonian generates a complete $C^1$ flow $\psi_t$.
We consider the Hamilton-Jacobi equation
\begin{equation}\tag{HJ}
H(x,du_x)=c,
\end{equation}
with a special emphasis on sub-solutions.
A function $u:M\lto \Rm$
is called a sub-solution of (HJ) if it is Lipschtiz and satisfies the 
inequality $H(x,du_x)\leq c$
at almost every point.
Note that this definition is equivalent to the 
notion of viscosity sub-solutions, see \cite{Fa:un}.
We denote by $C^{1,1}(M,\Rm)$ the set of differentiable functions
with Lipschitz differential.
The goal of the present paper is to present a short and 
selfcontained proof of:
\vs\\
\textbf{Theorem A}
\begin{itshape}
Let $H$ be a Tonelli Hamiltonian.
If the Hamilton-Jacobi equation (HJ)
has a sub-solution, then it has a $C^{1,1}$ sub-solution.
Moreover, the set of $C^{1,1}$ sub-solutions is dense
for the uniform topology in the set of sub-solutions.
\end{itshape}
\vs

Fathi and Siconolfi recently proved the existence of a $C^1$
sub-solution in \cite{FS:04}, see \cite{Massart} for the non-autonomous case.
Our result is optimal in the sense that examples
are known where a $C^{1,1}$ sub-solution exists, but no $C^2$ 
sub-solution, see Appendix A.
There exists a real number $\alpha(H)$, 
called the Ma\~n\'e critical value
in the literature, such
that the equation  (HJ) has sub-solutions if and only if
$c\geq \alpha(H)$.
One can prove 
the existence of smooth sub-solutions 
for $c>\alpha(H)$ by standard  regularization, see \cite{CIPP}.
As a consequence, our Theorem is relevant for the sub-solutions
of the critical equation $H(x,du_x)=\alpha(H)$, which are called the
critical sub-solutions of (HJ).
The study of the critical Hamilton-Jacobi equation
$H(x,du_x)=\alpha(H)$ is the core of Fathi's weak KAM theory.

A sub-solution $u$ is called strict on the open 
set $U\subset M$ if there exists a smooth non-negative 
function $V:M\lto \Rm$ which is positive on $U$ and such that $u$
is also a sub-solution of the equation 
$H(x,du_x)+V(x)=c$.
By applying the Theorem to the Hamiltonian $H+V$, we obtain:
\vs\\
\textbf{Addendum}
\begin{itshape}
If there exists a sub-solution of (HJ) which is strict
on the open set $U$, then there exists a $C^{1,1}$ sub-solution which
is strict on $U$.
\end{itshape}
\vs

We now expose some dynamical consequences of the main result,
which lead to a very short proof of the existence of
invariant sets contained in Lipschitz graphs:
\vs\\
\textbf{Theorem B}
\begin{itshape}
There exists a unique non-empty
 compact set $\tilde \mA(H)\subset T^*M$
with the following properties:
\begin{enumerate}
\item
$\tilde \mA(H)$ is  invariant  for the 
Hamiltonian flow, and 
$$\tilde \mA(H)\subset H^{-1}(\alpha(H)).$$
\item
For each $C^1$ critical sub-solution $u$ of (HJ),
we have 
$$\tilde \mA(H)\subset \Gamma_u:=\{(x,du_x) | x\in M\}.$$
\item
There exists a critical $C^{1,1}$ sub-solution $u$
of (HJ)  which is strict on the complement of
the projection $\mA(H)$ of $\tilde \mA(H)$ onto $M$.
\end{enumerate}
\end{itshape}
\vs

It is an easy consequence of Theorem B that the set $\tilde \mA(H)$
is a Lipschitz graph above $\mA(H)$ and is not empty.
We  explain in the course of the proof of Theorem B
in section \ref{Aubry} that
$\tilde \mA(H)$  is the set 
 usually called the Aubry set in the literature
(although it was introduced by John Mather). 

Let us quote explicitely the following:
\vs
\\
\textbf{Corollary}
\begin{itshape}
There exists a critical  $C^{1,1}$ sub-solution which is strict outside of
the projected Aubry set.
\end{itshape}
\vs

We give some examples in Appendix A, which explain
why $C^{1,1}$ regularity is optimal.
Theorem A is proved in Section 1, with the use of some
properties of semi-concave functions which are recalled
in Appendix B. Theorem B is proved  in Section 2.

I wish to thank Cedric Villani, whose questions
about  the geometry of optimal transportation
led me to the Proposition \ref{key} which is the key of the proof.

\section{The Lax-Oleinik semi-groups and sub-solutions}

We prove Theorem A.
It is necessary to start with more definitions.
We define the Lagrangian $L:TM\lto \Rm$ associated to $H$
by the relation
$$
L(x,v)=
\max_{p\in T_x^*M} p(v)-H(x,p).
$$
Then we define, for each $t\geq 0$, the function 
$A_t:M\times M\lto \Rm$ by
$$
A_t(x,y):=\min_{\gamma}
\int_0^t c+L(\gamma(s),\dot \gamma(s))ds
$$
where the minimum is taken on the set of curves
$\gamma\in C^2([0,t],M)$ which satisfy
$\gamma(0)=x$ and $\gamma(t)=y$.
Following Fathi, we define the Lax-Oleinik semi-groups
$T_t$ and $\breve T_t$ on $C^0(M,\Rm)$ by 
$$
T_t u(x)=\min_{y\in M}\big( u(y)+A_t(y,x)\big)
\;\;\text{ and }\;\;
\breve T_t u(x)=\max_{y\in M}\big( u(y)-A_t(x,y)\big).
$$
The following useful Lemma is proved in Fathi's book:

\begin{lem}
Given a Lipschitz function $u:M\lto \Rm$, the following properties are equivalent:
\begin{itemize}
\item $u$ is a sub-solution of (HJ).
\item The inequality
$u(y)-u(x)\leq A_t(x,y)
$
holds for each $t> 0$ and each $(x,y)\in M\times M$.
\item 
The function $[0,\infty[\ni t\lmto  T_tu(x)$ is non-decreasing for each $x\in M$.
\item
The function $[0,\infty[\ni t\lmto  \breve T_tu(x)$
 is non-increasing for each $x\in M$.
\end{itemize}
\end{lem}

An important consequence is that the semi-groups $T_t$
and $\breve T_t$ preserve the set of sub-solutions.
Another important property of these semigroups is that,
for each  $t>0$ and each  continuous function $u$,
the function $T_tu$ is semi-concave and the function
$\breve T_t u$ semi-convex,
see \cite{AMS,Fa:un} and  Appendix B for the definitions.
Recall  that a function is $C^{1,1}$
if and only if it is both semi-concave and semi-convex.

If $u$ is a sub-solution of (HJ),
then for each $s> 0$ and $t> 0$,
the function 
$T_s\breve T_t u$
is a sub-solution. We shall prove that, for each fixed 
$t>0$, this function is $C^{1,1}$ when $s$ is small enough.
Ludovic Rifford has pointed out to the author that this is a kind of
Lasry-Lions regularization, see \cite{LL}.
Since $\breve T_t u$ is semi-convex, Theorem A follows from
the following result, which may have other applications.

\begin{prop}\label{key}
Let $H$ be a Tonelli Hamiltonian.
For each semi-convex function $v$, the function 
$T_sv$ is $C^{1,1}$ for each  sufficiently small $s>0$. 
\end{prop}

\proof
In order to prove this proposition, it is enough to prove that
the function $T_s v$ is semi-convex for small $s$,
since we already know that it is semi-concave for all $s>0$.
This follows from two Lemmas:

\begin{lem}
For each  bounded  subset $F\subset C^2(M,\Rm)$
 there exists a 
time $s_0>0$ such that, for each $s\in [0,s_0]$,
the image $T_s(F)$ is a bounded subset of $C^2(M,\Rm)$ and
the following relation holds for 
all functions $f\in F$ and all $x\in M$
\begin{equation}\label{eq}
T_sf(x(s))=f(x)+\int_0^s c+L(x(t),\dot x(t)) dt,
\end{equation}
 where $x(t)$ is the curve 
$\pi \circ \psi_t(x,df(x))$
($\pi:T^*M\lto M$ is the projection and $\psi_t$
is the Hamiltonian flow).
\end{lem}
\proof
Let us  consider a $C^2$ function $f$ and the graph
$\Gamma_f\subset T^*M$ of its differential.
This graph is a $C^1$ Lagrangian manifold transversal to the
fibers.
It is known that, for $s\geq 0$ small enough, the Lagrangian manifold
$\psi_s(\Gamma_f)$ is the graph of a $C^2$
function, and that this $C^2$ function is $T_sf$.
Then, we have (\ref{eq}). 
The maximum time $s_0$ such that these properties hold
is uniform for families of functions
which are bounded in $C^2$ norm
(for then the associated graphs are contained
in a given compact set, and are uniformly transversal to
the verticals).
In addition, one can choose $s_0$ in such a way that
the set  $\{T_sf,s\in [0,s_0],f\in F\}$
is bounded in the $C^2$ topology,
(which amounts to say that the manifolds
$\psi_s(\Gamma_f)$ are uniformly transversal to the fibers).
 \qed

\begin{lem}
Let $v$ be a semi-convex function.
Then there exists a bounded subset $F\subset C^{2}(M,\Rm)$
and a time $s_0>0$ such that 
$$
T_sv=\max _{f\in F} T_sf
$$
for all $s\in [0,s_0]$, hence $T_sv$ is a semi-convex function
for $s\geq 0$ small enough.\vs\\
\end{lem}
\proof
If $v$ is semi-convex, then there exists a 
bounded subset $F\subset C^2(M,\Rm)$
such that
$v=\max _F f$ and such that
 for each $x$ and each $p\in \partial^-v(x)$ 
 (the set of proximal sub-differentials of v at point $x$, 
 see Appendix B),
there exists a function $f\in F$ satisfying
$(f(x),df(x))=(v(x),p)$, see Appendix B.
Let us fix from now on such a family $F$ of functions,
and consider the time $s_0$ associated to this family
by the first Lemma.
Notice that 
$$
T_sv\geq \sup _{f\in F} T_sf
$$
for all $s$, because for each $f\in F$ we have $f\leq v$ hence
$T_sf\leq T_sv$.
In order to prove that the equality holds
for $s\in [0,s_0]$, 
let us fix a point $x\in M$.
There exists a point $y$ such that 
$$
T_sv(x)=v(y)+A_s(y,x).
$$
Now let $(x(t),p(t)):[0,s]\lto T^*M$
be a Hamiltonian trajectory which is optimal 
for $A_s(y,x)$.
We mean that $x(0)=y$, $x(s)=x$,
and 
$$
A_s(y,x)=\int_0^s c+L(x(t),\dot x(t)) dt.
$$
It is known (see \cite{Fa:un,AMS}) that
$-p(0)$ is then a proximal super-differential of the function
$z\lmto A_s(z,x)$ at point $y$.
 Since the function 
$z\lmto u(z)+A_s(z,x)$  is minimal at $y$,
the  linear form $p(0)$ is
 a proximal  sub-differential of the function $u$ at point 
$y$.
Let us consider a function $f\in F$
such that $(f(y),df(y))=(u(y),p(0))$.
Then we have  $(x(t),p(t))=\psi_t(y,df(y))$
and, by the first Lemma, 
$$
T_sf(x)=T_sf(x(s))=f(y)+\int_0^sc+ L(x(t),\dot x(t)) dt
=u(y)+A_s(y,x)=T_su(x).
$$
We have proved that, for each point $x\in M$,
there exists a function $f\in F$ such that 
$T_sf(x)=T_su(x)$.
This ends the proof.
\qed

The proof also implies:
\begin{cor}\label{cor}
If $u$ is a $C^{1,1}$ sub-solution, then there exists $\epsilon>0$
such that $T_{t}u$ and $\breve T_{t}u$
are $C^{1,1}$ sub-solutions when $t\in [0,\epsilon]$. In addition, we have, for these values of $t$,
$$
\Gamma_u=\psi_{t}\big(\Gamma_{\breve T_{t}u}\big)
=\psi_{-t}\big(\Gamma_{ T_{t}u}\big)
$$
where $\Gamma_f$ is the graph of the differential of $f$.
\end{cor}

\section{The Aubry set}\label{Aubry}
In this section, we consider only the critical case $c=\alpha(H)$,
and prove Theorem B.
Let us first define the projected Aubry set $\mA(H)\subset M$.
This is the set of points $x\in M$ such that
$H(x,du_x)=\alpha(H)$ for each $C^1$ sub-solution $u$.
A similar definition is given in \cite{FS:05}.

 \begin{lem}
If $u_1$ and $u_2$ are two critial $C^1$ sub-solutions,
then $du_1=du_2$ on $\mA(H)$.
 \end{lem}
 \proof
If $du_1(x)\neq du_2(x)$,
then, by the strict convexity of $H$, the function
$(u_1+u_2)/2$ is a $C^1$ critical sub-solution which is 
strict at $x$. 
This implies that $x$ does not belong to $\mA(H)$.
\qed

As a consequence, we can define in a natural way the set
$$
\tilde \mA(H):=
\{(x,du_x)| x\in \mA(H)\}
$$
where $u$ is any $C^1$ critical sub-solution.

\begin{lem}
There exists a $C^{1,1}$ critical sub-solution 
$u$ which is strict outside of $\mA(H)$.
\end{lem}

\proof
By the Addendum of Theorem A, it is enough to prove 
that there exists a  critical sub-solution
which is strict outside of $\mA(H)$.
Since $C^1(M,\Rm)$ is separable, the set of  
critical $C^1$ sub-solutions
of (HJ) endowed with the $C^1$   norm is separable.
As a consequence, there exists a dense  sequence $u_n$
 of $C^1$ critical sub-solutions.
 The $C^1$ function 
$$
u(x):=\sum_{n=1}^{\infty} \frac{u_n(x)}{2^n}
$$
is a $C^1$ critical sub-solution of (HJ)
which is strict outside of $\mA(H)$.
Indeed, for each point $x\not \in \mA(H)$,
there exists a 
$C^1$ critical sub-solution $v$ such that
$H(x,dv_x)<\alpha(H)$.
Since the sequence $u_n$ is dense for the $C^1$
topology, we conclude that $H(x,du_n(x))<\alpha(H)$
for some $n$. The desired conclusion follows by the 
convexity of $H$.
\qed

This Proposition implies that $\mA(H)$ is not empty.
Otherwise, there would exist a critical sub-solution strict on $M$,
 which is a contradiction.

\begin{prop}
The set $\tilde \mA(H)$ is invariant.
\end{prop}

\proof
Let us choose a $C^{1,1}$ critical sub-solution $u$
which is strict outside of the projected Aubry set.
We have $\tilde \mA(H)=\Gamma_u\cap H^{-1}(\alpha(H))$.
Let $\epsilon$ be given by corollary \ref{cor}.
We claim that $\psi_t(\tilde \mA(H))=\tilde \mA(H)$
for all $t\in [-\epsilon,\epsilon]$,
where $\psi_t$ is the Hamiltonian flow. This claim clearly implies 
the desired result.
Let $(x,du_x)$ be a point of $\tilde \mA(H)$ and
$t\in [-\epsilon,\epsilon]$. 
Let us denote by $(y,dv_y)$ the point  $\psi_t(x,du_x)$,
where $v:= T_tu$.
Since $v$ is a critical sub-solution,
we  have  $(y,dv_y)\in \tilde \mA(H)$
provided $y\in \mA(H)$.
In order to prove this inclusion, 
we denote by $w$ the function 
$w:=\breve T_tu$, which is a $C^{1,1}$ critical sub-solution.
Since $x\in \mA(H)$, we have $du_x=dw_x$.
This implies that $\psi_t(x,dw_x)=\psi_t(x,du_x)=(y,dv_y)$.
Since $\psi_t(\Gamma_w)=\Gamma_u$,
this implies that $dv_y=du_y$, and, by energy conservation,
that $H(y,du_y)=\alpha(H)$. 
Since the sub-solution $u$ is strict outside
of $\mA(H)$, we conclude that $y\in \mA(H)$.
\qed
The usual definition of the Aubry set is based on the notion of calibrated curve.
\footnote{here start the modifications compared to the published version} 
The curve $x(t)$ is said calibrated by the sub-solution $u$
if the equality
\begin{equation}\label{calibration}
u(x(t))-u(x(s))=\int_s^t \alpha(H)+L(x(\sigma),\dot x(\sigma))d\sigma
\end{equation}
holds for each $s\leq t$.
If $x(t)$ is calibrated by $u$, then it is an extremal, which means that 
there exists a trajectory $(x(t),p(t))$ of the Hamiltonian system above the curve $x(t)$.
In the book of Fathi \cite{Fa:un}, the Aubry set is defined 
as the union, in $T^*M$, of all the images of  the Hamiltonian orbits
$(x(t),p(t)):\Rm \lto T^*M$ which are calibrated by all subsolutions.
This definition is equivalent to the ones previously given by John Mather.
Denoting temporarily  by $\tilde \mF(H)$ this set defined by Fathi,
we have 
 $$\tilde \mA(H)= \tilde \mF(H),
$$
which explains our terminology.
The inclusion $\tilde \mA(H)\subset \tilde \mF(H)$
follows from:

\textbf {Claim.} 
\begin{itshape}
Let $(x(t),p(t))$ be a Hamiltonian orbit contained in the
invariant set $\tilde \mA(H)$.
Then the curve $x(t)$ is calibrated by all the critical sub-solutions.
\end{itshape}

\proof
Since the $C^{1,1}$ critical sub-solutions are dense in the set of critical sub-solutions for the uniform topology,
it is sufficient to prove the statement for $C^{1,1}$ critical sub-solutions.
Let $u$ be such a sub-solution.
The point $(x(t),p(t))$ belongs to $\tilde \mA(H)$, hence it belongs
to the graph of $du$, and therefore $p(t)=du_{x(t)}$.
We thus  have
$$
L(x(t),\dot x(t))=p(t)\cdot \dot x(t)-H(x(t), p(t))=du_{x(t)}\cdot \dot x(t)-\alpha(H)
$$
and the desired formula follows by integration.
\qed
Conversely, in order to prove the inclusion $\tilde \mF(H)\subset \tilde \mA(H)$, we neeed:

\textbf {Claim.} 
\begin{itshape}
Let 
 $(x(t),p(t))$ be a Hamiltonian orbit contained in $\tilde \mF(H)$, then,
for each $C^1$ sub-solution $u$, we have
 $p(t)=du_{x(t)}$ and $H(x(t),p(t))=\alpha(H)$ for each $t\in \Rm$.
\end{itshape}

\proof
It is a  classical fact of Weak KAM theory, proved for example in \cite{Fa:un}
that $p(t))=du_{x(t)}$ when $(x(t),p(t))$ is an orbit calibrated by the $C^1$ sub-solution 
$u$ (actually this equality holds even if $u$ is not $C^1$).
We conclude that all $C^1$ sub-solutions satisfy
$du_{x(t)}=p(t)$.
Assume that $H(x(t),p(t))<\alpha(H)$. Then it is possible to perturb slightly the function $u$
by a function $v$ which is still a critical sub-solution and satisfies $dv_{x(0)}\neq p(0)$, a contradiction.
\qed

\begin{prop}
If $u$ is a critical sub-solution (not necessarily $C^1$), then
$T_tu(x)=\breve T_tu(x)=u(x)$ for all $t\geq 0$ and $x\in \mA(H)$.
Therefore, if $u$ is a critical sub-solution, there 
exists a $C^{1,1}$ sub-solution which coincides with $u$ 
on $\mA(H)$. 
\end{prop}

\proof
The second part of the statement clearly follows from the first
part: just take $T_{\epsilon}\breve T_t u$, which is equal to $u$
on $\mA(H)$.
So we have to prove the first part of the statement.
Once again, it is enough to prove it when $u$ is $C^{1,1}$.
We now make this additional assumption.
Let $x$ be a point in $\mA(H)$, and let $(x(t), p(t))$ be the  orbit of
the point $(x, du_x)$.
Since the curve $x(t)$ is calibrated by $u$, we have 
$$
u(x)=u(x(-t))+\int_{-t}^0 \alpha(H)+L(x(s),\dot x(s))ds\geq T^t u(x)
$$
for all $t\geq 0$. Since $u$ is  sub-solution, we conclude that $T^tu(x)=u(x)$.
The proof concerning $\breve T$ is similar.
\qed

\appendix
\section{Examples}\label{examples}
\subsection{Mechanical Hamiltonian system}
Let us consider the case
$$
H(x,p)=\frac{1}{2}\|p\|_x^2+V(x)
$$
where $\|.\|_x$ is a Riemaniann metric on $M$
and $V$ is a smooth function on $M$.
Then it is easy to see that 
$\alpha(H)=\max V$,
and that there exists a smooth sub-solution to (HJ):
any constant function is such a sub-solution!

\subsection{Non-existence of a $C^2$ sub-solution}
Let us now specialise to $M=\Tm$,
and consider the Hamiltonian 
$$
H_P(x,p)=\frac{1}{2}(p+P)^2-\sin^2 (\pi x)
$$
depending on the real parameter $P$.
For $P=0$ this is a Mechanical system as discussed above,
and the constants are sub-solutions of (HJ).
Let $X(x):\Tm\lto \Rm$ be the function
such that $X(x)=\sin(\pi x)$ for $x\in [0,1]$.
Let us set
$$
a=\frac{2}{\pi}
=\int_{\Tm} X(x) dx.
$$
The reader can check easily that 
$\alpha(H_P)=0$ for $P\in [-a,a]$.
 For each $P\in ]-a,a[$,
the equation (HJ) has smooth sub-solutions.
For these values of $P$, the Aubry set is the fixed   point $(0,-P)$.

However, for $P=a$, there is one and only one 
critical sub-solution
of (HJ), which turns out to be a solution.
It is given by the primitive of the function 
$X-a$.
This function is $C^{1,1}$  but not $C^2$.
Note that the Aubry set, then, is not reduced to the 
hyperbolic fixed point $(0,-a)$ but is the whole graph of $X-a$.

\section{Semi-concave functions}\label{semiconcave}
We recall some useful facts on semi-concave functions, 
see for example \cite{CaSi,Fa:un}
for more material.
In all this section, $M$ is a compact manifold of dimension $d$.
It is useful to fix once and for all 
 a finite atlas
$\Phi$ of $M$ composed of
charts  $\varphi:B_3\lto M$, where 
$B_r$ is the open  ball of radius $r$
centered at zero in $\Rm^d$.
We assume that the sets $\varphi(B_1),\varphi\in \Phi$ cover $M$.
A family $F$  of  $C^2$ functions
is said bounded if there exists a constant $C>0$
such that
$$
\|d^2(u\circ \varphi)_x\|\leq C
$$
for all $x\in B_1, \varphi\in \Phi, u\in F$.
Note that  a bounded family is not required
to be bounded in $C^0$ norm, but will automatically
be bounded in $C^1$ norm and thus equi-Lipschitz.
The notion of bounded family of functions
does not depend on the atlas $\Phi$.

A function $u:M\lto \Rm$ is called semi-concave 
if there exists a bounded subset $F_u$ of the set
$C^2(M,\Rm)$  such that
$$
u=\inf_{f\in F_u}f.
$$
A semi-concave function is Lipschitz. 
We say that the linear form $p\in T_xM$ is a
proximal super-differential
of the function $u$ at point $x$ if there exists a $C^2$
function $f$ 
such that $f-u$ has a minimum at $x$ and $df_x=p$.
We denote by $\partial^+u(x)$ the set of proximal superdifferentials 
of $u$ at $x$.
We say that a linear form $p\in T_xM$
is a $K$-super-differential of the function
$u$ at point $x$ if for each chart
$\varphi\in \Phi$ and each  point $y\in B_2$
satisfying
$\varphi(y)=x$, the inequality
$$
u\circ\varphi(z)-u\circ\varphi(y)
\leq p\circ d\varphi_y(z-y)+K\|z-y\|^2
$$
holds for each $z\in B_2$.
A function $u$ on $M$ is called
$K$-semi-concave if it
has a $K$-super-differential at each point.
It is equivalent
to require that, for each $\varphi\in \Phi$,
the function
$$
u\circ \varphi(y)-K\|y\|^2
$$
is concave on $B_2$.
As a consequence, if $u$ is $K$-semi-concave
and if $p$ is a proximal super-differential of $u$
at $x$, then $p$ is a $K$-super-differential of
$u$ at $x$.

\begin{prop}
A function $u$ is semi-concave if and only if
there exists a number  $K>0$ such that 
$u$ is $K$-semi-concave.
Then,  there  exists a bounded subset $F\subset C^2(M,\Rm)$
 such that 
$$
u=\min _{f\in F} f
$$
and, for each point $x\in M$ and each super-differential
$p$ of $u$ at $x$, there exists a function $f\in F$
such that $(f(x),df(x))=(u(x),p)$.
\end{prop}

\proof
Let us consider a smooth function 
$g:\Rm^d\lto \Rm$ such that 
$0\leq g \leq 1$, and such that
$g=0$ outside of $B_2$ and $g=1$ inside $B_1$.
Let us associate, to each chart $\varphi\in \Phi$,
and each point $(x,p)\in T_xM$
satisfying $x\in \varphi(B_1)$,
 the function $f_{x,p,\varphi}:M\lto \Rm$
defined by
$$
f_{x,p,\varphi}\circ \varphi(z)
:=g(z)\big(u(x)+p\circ d\varphi_y(z-y)+K\|z-y\|^2\big)
+(1-g(z))\max u
$$
for $z\in B_2$, where $y=\varphi^{-1}(x)$,
and 
$f_{x,p,\varphi}=\max u$
outside of  $\varphi(B_2)$.
The functions $f_{x,p,\varphi}, \varphi\in \Phi,
x\in \varphi(B_1), p\in \partial^+u(x)$
form a bounded subset $F$ of $C^2(M,\Rm)$.
It is easy to check that 
$f=\min _{f\in F}f$.
\qed

A function $u$ is called semi-convex if $-u$ is semi-concave.
\begin{prop}
A function is $C^{1,1}$ if and only if it is semi-concave and
semi-convex.
\end{prop}

A very elementary proof of this statement is given in the
book of Fathi. Another proof is given in \cite{CaSi},
Corollary 3.3.8.

\small
\bibliographystyle{amsplain}
\providecommand{\bysame}{\leavevmode\hbox
to3em{\hrulefill}\thinspace}

\end{document}